\documentclass[11pt]{amsart}
\usepackage{epsfig}
\usepackage{amssymb}
\usepackage{amscd}
\usepackage{epsf}

\long\def\comment#1\endcomment{\relax}

\newcounter{subsubsubsection}

\makeatletter
\@addtoreset{subsubsubsection}{subsubsection}

\newcommand{\subsubsubsection}[1]{\par\addtocounter{subsubsubsection}{1}\smallskip
\noindent{\thesubsubsection.\arabic{subsubsubsection}.\bfseries\ #1.}}

\newcommand{\subsubsubsectionb}[1]{\par\addtocounter{subsubsubsection}{1}\smallskip
\noindent{\thesubsubsection.\arabic{subsubsubsection}.\bfseries\ #1}}

\makeatother

\newcommand{\sevafig}[4]{\begin{figure}[h]\centerline{
 \epsfig{file=td#1.eps,width=#2,angle=#3}}
\vskip-1cm
\caption{#4}\end{figure}}

\newcommand{\one}{{\mathbf1}}

\newcommand{\U}{\mathcal U}

\newcommand{\PBW}{{PBW}}

\newcommand{\simto}{\overset{\sim}{\to}}

\newcommand{\ndot}{\bullet}

\newcommand{\g}{{\mathfrak g}}

\newcommand{\gl}{{\mathfrak{gl}}}

\def\matho#1{\mathop{\mathrm{#1}}}

\DeclareMathSizes{11.1}{10}{8}{6}

\newcommand{\opp}{{\mathrm{opp}}}

\newcommand{\Tr}{{\matho{Tr}}}

\newcommand{\Id}{\matho{Id}\nolimits}
\newcommand{\ad}{\matho{ad}\nolimits}
\newcommand{\Log}{\matho{Log}}

\newcommand{\strange}{{\mathrm{strange}}}
\newcommand{\poly}{{\mathrm{poly}}}

\newcommand{\Hoch}{\mathrm{Hoch}}

\newcommand{\C}{\mathbb C}
\newcommand{\R}{\mathbb R}
\newcommand{\Z}{\mathbb Z}
\newcommand{\D}{\mathcal D}

\newtheorem*{theorem}{Theorem}

\newtheorem*{corollary}{Corollary}
\newtheorem*{conjecture}{Conjecture}

\theoremstyle{remark}
\newtheorem*{remark}{Remark}

\theoremstyle{definition}

\author{Boris Shoikhet}
\title
{Tsygan formality and Duflo formula}
\date{2000}
\address{FIM, ETH-Zentrum, CH-8092 Z\"urich, SWITZERLAND}
\email{borya@mccme.ru}

\begin{document}

\sloppy

\maketitle

\begin{abstract}
We prove the 0-(co)homology part of the conjecture on the cup-products on tangent cohomology in the Tsygan formality [Sh2]. We discuss its applications to the Duflo formula.
\end{abstract}

\section*{A short introduction}

The Tsygan formality conjecture for chains~[Ts] was proven in
the author's work~[Sh2] by an explicit construction of suitable
Kontsevich-type integrals. This paper is a further development of
ideas of~[Sh2]. We will freely use the notations and results
of~[Sh2]. In~[Sh2] we formulated a conjecture on ``the cup-products
on tangent cohomology'', which is a version of the analogous
Kontsevich theorem from Section~8 of~[K]. Here we prove this conjecture for 0-(co)homology.

\section{The classical Duflo formula and the generalized Duflo
formula}

\subsection{}
Let $\g$ be a finite-dimensional Lie algebra, $S^\ndot(\g)$ and
$U(\g)$ be its symmetric and universal enveloping algebra. They
are not isomorphic as algebras $S^\ndot(\g)$ is a commutative
algebra and $U(\g)$ is a non-commutative algebra. We can
consider both spaces $S^\ndot(\g)$ and $U(\g)$ as $\g$-modules
with the adjoint action for $S^\ndot(\g)$ and the action
$g\cdot\omega=g\otimes\omega-\omega\otimes g$ for~$U(\g)$ (here
$g\in\g$ and $\omega\in U(\g)$). It is clear that these
$g$-modules are isomorphic, the isomorphism is the classical
Poincar\'e--Birkhoff--Witt map:
\begin{equation}
\varphi_\PBW(g_1\cdot\dots\cdot g_k)=\frac1{k!}
\sum_{\sigma\in\Sigma_k}g_{\sigma(1)}\otimes\dots\otimes g_{\sigma(k)}
\end{equation}
($g_1,\dots,g_k\in\g$).

The Duflo theorem~[D] states that the invariants
$[S^\ndot(\g)]^\g$ and $[U(\g)]^\g$ are isomorphic \emph{as
algebras}. The Duflo formula is a canonical formula for this
isomorphism. We recall it here.

For any $k\ge1$, there exists a canonical element in
$[S^k(\g^*)]^\g$. It is the symmetrization of the map
$$
g\mapsto\Tr|_\g\ad^kg\quad(g\in\g).
$$
We denote this element in $[S^k(\g^*)]^\g$ by~$\Tr_k$.  We can consider an
element from $S^k(\g^*)$ as a differential operator of the
$k$-th order with constant coefficients, acting
on~$S^\ndot(\g)$. (Thus, an element from~$\g^*$ is a derivation
of~$S^\ndot(\g)$). (It was a
conjecture of M.~Duflo that the operators corresponding to
$S^k(\g^*)^\g$ are \emph{zero} for
odd~$k$ and any finite-dimensional Lie algebra~$\g$; this
conjecture was proven recently in [AB]).

Define the map $\varphi_\strange\colon S^\ndot(\g)\to
S^\ndot(\g)$ by the formula:
\begin{equation}
\varphi_\strange=\exp\left(\sum_{k\ge1}\alpha_{2k}\cdot\Tr_{2k}\right)
\end{equation}
where
\begin{equation}
\sum_{k\ge1}\alpha_{2k}\cdot x^{2k}=\frac12\Log\frac{e^{\frac
x2}-e^{-\frac x2}}x.
\end{equation}
The map $\varphi_\strange$ is well-defined on~$S^\ndot(\g)$ (in
the sense that we have no problems with divergences), because
$\Tr_{2k}(\omega)\equiv0$ for a fixed $\omega\in S^\ndot(\g)$
and for a sufficiently large~$k$. The map $\varphi_\strange$ is
a map of $\g$-modules, because the operators $\Tr_{2k}$ are
invariant.

\begin{theorem}[M.~Duflo, {[D]}]
The restriction of the map
$\varphi_D=\varphi_\PBW\circ\varphi_\strange$ to the invariants
$[S^\ndot(\g)]^\g$ defines a map of \emph{algebras}
$\varphi_D\colon[S^\ndot(\g)]^\g\to[U(\g)]^\g$.
\end{theorem}

\subsubsection{}
M.~Kontsevich deduced from his theorem on cup-products on the
tangent cohomology~[K] the following generalization of the Duflo
theorem.

\begin{theorem}
There exists a canonical map $\tilde{\varphi_D}\colon
H^\ndot(\g;S^\ndot(\g))\to H^\ndot(\g;U(\g))$
which is a map of associative algebras. Its restriction to
$H^0(\g;S^\ndot(\g))$ coincides with the Duflo map $\varphi_D$.
This result holds also for any $\Z$-graded
finite-dimensional Lie algebra~$\g$.
\end{theorem}

Recall, that $M^\g=H^0(\g;M)$ for any $\g$-module~$M$.

This map $\tilde{\varphi_D}$ is given as the tangent map to the Kontsevich
$L_\infty$ formality morphism at the solution to the Maurer-Cartan equation
corresponding to the Kostant-Kirillov Poisson structure on $\g^*$.
In the case of $H^0(\g;S^\ndot(\g))$ this tangent map can be computed (not so
easy, by comparing with the Duflo formula for $gl_N$ in the Kontsevich
original approach). In this case all the graphs are unions of so called wheels.
For higher cohomology $H^k(\g;S^\ndot(\g))$, $k\ge 1$, many other graphs
besides the wheels appear, and it seems that any computation of the
Feynmann(Kontsevich) weights of these other graphs is impossible.
Nevertheless, we can prove the following result.

\subsubsubsection{}
\begin{theorem}
Denote by $\varphi_D^\ndot\colon H^\ndot(\g;S^\ndot(\g))\to H^\ndot(\g;U(\g))$
the map induced by the map of $\g$-modules $\varphi_D\colon S^\ndot(\g)\to
U(\g)$. Then the map $\varphi_D^\ndot$ is a map (an isomorphism) of
associative (graded commutative) algebras.
\begin{proof}
We just sketch the proof here. The complete proof will appear somewhere.
This proof is based on an unpublished joint paper with Maxim Kontsevich.

Consider the space $V=\g[1]$. It is a $\mathbb{Z}$-graded vector
space. Consider the $L_\infty$ formality morphism on it. The polyvector
fields $T_{\poly}(V)$ is isomorphic to $T_\poly (\g^*)\simeq
\wedge^\ndot(\g^*)\otimes S^\ndot(\g)$. (This phenomenon can be considered
roughly as a kind of the Koszul duality). There is an odd vector field $Q$
on $\g[1]$ such that $Q^2=0$ (which generates the cochain differential). In
coordinates, $Q=\sum_{i,j,k=1}^{\dim\g}c_{ij}^k\xi_i \xi_j
\frac{\partial}{\partial \xi_k}$ where $c_{ij}^k$ are the structure
constants of the Lie algebra $\g$ in some basis $x_i$ and $\xi_i$ are the
odd coordinates on $\g[1]$ corresponding to $x_i$.

We want know to localize the formality morphism on $V$ at the solution to
the Maurer-Cartan equation $Q\in [T_\poly(V)]^1$.
We claim that the only graphs which appear are unions of the wheels.


It follows
from [K], Lemma7.3.3.1(1). Note that here the wheels are not the same wheels
as for $\g^*$: here we have one outgoing edge and two incoming edges for
each vertex whence for $\g^*$ we have two outgoing edges and one incoming. It
reflects the fact that $Q$ is a (quadratic) vector field whence the
Kostant-Kirillov Poisson structure $\alpha=\sum_{i,j,k=1}^{\dim \g}c_{ij}^k
x_k\frac{\partial}{\partial x_i}\wedge\frac{\partial}{\partial}{\partial
x_j}$ is a (linear) bivector field.

It is not straightforward to compute the Kontsevich weights of these wheels
corresponding to $\g[1]$ but it turns out it is possible. The answer is
exactly the formula in the Theorem.

\end{proof}
\end{theorem}

Later in this paper we consider only the case of 0-th cohomology, to
simplify the exposition. While the Kontsevich claim on
the cup-products was proved recently for higher cohomology in [MT], in
the case of Tsygan formality we prove the corresponding theorem on
cup-products for 0-cohomology only. It seems, however, that the technique
developed in [MT] can be used in this situation.

\subsection{}

In [Sh2] we proposed the following conjecture:

\begin{conjecture}
Denote by $\varphi_{D\ndot}\colon
H_\ndot(\g;S^\ndot(\g))\to H_\ndot(\g;U(\g))$ the map defined by the map of
$\g$-modules
$\varphi_D\colon S^\ndot(\g)\to U(g)$.
Then the map $\varphi_{D\ndot}\colon
H_\ndot(\g;S^\ndot(\g))\to H_\ndot(\g;U(\g))$ is a
map of modules from the $H^\ndot(\g;S^\ndot(\g))^\opp$-module
$H_\ndot(\g;S(\g))$ to the $(H^\ndot(\g;U(\g))^\opp$-module
$H_\ndot(\g;U(\g))$. It means that for any $\alpha\in
H^\ndot(\g;S(\g))$ and any $\beta\in H^\ndot(\g;S(\g))$ one
has\emph:
$$
\varphi_{D\ndot}(\alpha\blacklozenge\beta)=\varphi_{D}^\ndot
(\alpha)\bigstar\varphi_{D\ndot}(\beta).
$$
\end{conjecture}
Here we denote by $\blacklozenge$ the canonical action
of $\alpha\in H^\ndot(\g;S^\ndot(\g))$ on $\beta\in H_\ndot(\g;S^\ndot(\g))$
and by $\bigstar$ the action of $H^\ndot(\g,U(\g))$ on $H_\ndot(\g; U(g))$.
As in general, cohomology forms an algebra, and homology forms a module over
it.

Recall, that $H_0(\g;M)=M/\g M=M_\g$ is the space of
coinvariants. For $0$-cohomology this conjecture states that
$(S^\ndot(\g))^\g$-module $(S^\ndot(\g))_\g$ and
$(U(\g))^\g$-module $(U(\g))_\g$ are isomorphic by means of the
Duflo map~$\varphi_D$.

\subsubsection{}

We prove here the following statement:

\begin{theorem}
For any finite-dimensional Lie algebra~$\g$ \emph(or any
finite-dimensional $\Z$-graded Lie algebra~$\g$\emph) one
has\emph:
$$
\varphi_D(\alpha\cdot\beta+c(\alpha,\beta))=\varphi_D(\alpha)\star\varphi_D(\beta)
$$
where $\alpha\in[S^\ndot(\g)]^\g$, $\beta\in S^\ndot(\g)$,
$*$ is the product in~$U(\g)$, and $c(\alpha,\beta)\in \{\g,S^\ndot(\g)\}$.
\end{theorem}

As well, we obtain an explicit formula for $c(\alpha,\beta)$.

It is clear that this theorem implies Conjecture above for 0-(co)homology.
For a semisimple Lie algebra $\g$, this theorem is equivalent to the Duflo formula because of the decompositions
\begin{align}
S^\ndot(\g)&=[S^\ndot(\g)]^\g\oplus\{S(\g),S(\g)\},\\
U(\g)&=[U(\g)]^\g\oplus[U(\g),U(\g)].
\end{align}
which hold for any semisimple Lie algebra $\g$.

For an arbitrary Lie algebra $\g$, this theorem is not a corollary of the Duflo formula, and it is a new fact about the Duflo map.

\section{The theorem on cup-products in Tsygan formality}

Here we prove the conjecture on the cup-products in the Tsygan
formality~[Sh2] for $0$-cohomology. We use the notations
from~[Sh2].

This conjecture is analogous to the Kontsevich theorem on a
cup-products in~[K], Section~8. It would be helpful for reader
to know the Kontsevich's proof. It is proven for 0-th tangent cohomology in [K], and in
[MT] in the general case.

\subsection{}
Recall that the Kontsevich $L_\infty$-morphism
$\U\colon T^\ndot_\poly(\R^d)\to\D^\ndot_\poly(\R^d)$ (see~[K])
and the Lie derivatives $L_\Psi\colon C_\ndot(A,A)\to
C_\ndot(A,A)$ allows to define a
$T^\ndot_\poly(\R^d)$-$L_\infty$-module structure on the chain Hochschild complex
$C_\ndot(A,A)$, $A=C^\infty(\R^d)$ (see for
details~[T], [Sh2],Section~1). Thus, we have two
$L_\infty$-modules over $T^\ndot_\poly(\R^d)$: these are
$C_\ndot(A,A)$ and~$\Omega^\ndot(\R^d)$, the differential forms
on $\R^d$ with zero differential and usual module structure
over~$T^\ndot_\poly$, defined with through the Lie derivatives
$L_\gamma=i_\gamma\circ d\pm d\circ i_\gamma$ (see~[T]).

In [Sh2] we constructed an $L_\infty$-morphism of
$L_\infty$-modules over~$T^\ndot_\poly(\R^d)$, $\hat\U\colon
C_\ndot(A,A)\to\Omega^\ndot(\R^d)$. Its Taylor components are maps
$$
\hat\U_k\colon\Lambda^kT^\ndot_\poly(\R^d)\otimes
C_\ndot(A,A)\to\Omega^\ndot(\R^d)[-k].
$$
They are constructed as sums over all admissible graphs, \dots\  etc.

Then for any solution $\pi$ of the Maurer--Cartan equation
in~$T^\ndot_\poly(\R^d)$, i.e.\ of the equation $[\pi,\pi]=0$,
one can define the \emph{tangent map}
$$
T_\pi\hat\U\colon T_\pi C_\ndot(A,A)\to T_\pi\Omega^\ndot(\R^d)
$$
where $T_\pi C_\ndot(A,A)=C_\ndot(A_*,A_*)$ with $A_*$ the
Kontsevich deformation quantization, and
$T_\pi\Omega^\ndot(\R^d)=\{\Omega^\ndot(\R^d),L_\pi\}$
(see~[Sh2], Section~3 for details). The main property of the
map~$T_\pi\hat\U$, which follows immediately from the
$L_\infty$-morphism equations, is that $T_\pi\hat\U$ is a
\emph{map of complexes}. In degree~$0$, $T^0_\pi
C_\ndot(A,A)=A_*$ (considered as a vector space), and
$T^0_\pi\Omega^\ndot(\R^d)=A$. In degree~$0$ we obtain a map
$T_\pi\hat\U\colon A_*\simto A$ which is a map of homology,
i.e.\ $T_\pi\hat\U$ induces a map $A_*/[A_*,A_*]\simto
A/\{A,A\}$ (here $[A_*,A_*]$ is the commutant of the deformed
algebra, and $\{A,A\}$ is the commutant with respect to the
Poisson bracket). See~[Sh2], Section~3 for details.

The last property means that $T_\pi\hat\U$ maps $[A_*,A_*]$ to
$\{A,A\}$.

Now we are going to prove the following result.

\begin{theorem}
For any Poisson structure~$\pi$, any function $\alpha\in A$ such
that $[\alpha,\pi]=0$ and any $\beta\in A_*$ one has\emph:
\begin{equation}
T_\pi\hat\U(((T_\pi\U)\alpha)*\beta)=\alpha \cdot T_\pi\hat\U(\beta)+c(\alpha,\beta).
\end{equation}
Here $T_\pi\U$ is the tangent map with respect to the Kontsevich morphism,
$$
T_\pi\U\colon\bigl\{T^\ndot_\poly[1],d=\ad\pi\bigr\}\to C^\ndot(A_*,A_*),
$$
$*$ in the l.h.s. of \emph{(6)} is the Kontsevich
star-product, see~\emph{[K]}, Section~\emph8, and $c(\alpha,\beta)\in \{A,A\}$, the Poisson commutant of the algebra $A$.
\end{theorem}

\subsection{}
First of all, recall the definitions of the tangent
maps~$T_\pi\U$,~$T_\pi\hat\U$. The case of the Kontsevich
formality (i.e.\ the case of $L_\infty$-morphism between dg Lie
algebras) is simpler. We have:
\begin{equation}
(T_\pi\U)(x)=\U_1(x)+\U_2(x,\pi)+\frac12\U_3(x,\pi,\pi)+
\frac1{3!}\U_4(x,\pi,\pi,\pi)+\dots
\end{equation}
It is a map of complexes
$$
T_\pi\U\colon\bigl\{T^\ndot_\poly(\R^d)[1],d=\ad\pi\bigr\}\to
\bigl\{\D^\ndot_\poly(\R^d)[1],d=d_\Hoch+\ad\tilde\pi\bigr\}
$$
where
$$
\tilde\pi=\U_1(\pi)+\frac12\U_2(\pi,\pi)+\frac1{3!}\U_3(\pi,\pi,\pi)+\dots
$$
is the Kontsevich solution of the Maurer--Cartan equation
in~$\D^\ndot_\poly(\R^d)$. The last complex can be identified
with the
complex~$C^\ndot(A_*,A_*)$. More precisely, we should set
$\pi:=\hbar\pi$, where $\hbar$ is a formal parameter.

In the case of the Tsygan formality
\begin{equation}
(T_\pi\hat\U)(\omega)=\hat\U_1(\omega)+\hat\U_2(\pi,\omega)+
\frac12\hat\U_3(\pi,\pi,\omega)+\dots
\end{equation}
It is a map of the complexes
$$
T_\pi\hat\U\colon T_\pi C_\ndot(A,A)\to T_\pi\Omega^\ndot(\R^d)
$$
where for an $L_\infty$-module $M$ over dg Lie
algebra~$\g^\ndot$, and a solution~$\pi$ of the Maurer--Cartan
equation in~$\g^\ndot$, the differential in~$T_\pi M$ is equal to
\begin{equation}
d\omega=\phi_0(\omega)+\phi_1(\pi,\omega)+\frac12\phi_2(\pi,\pi,\omega)+\dots
\end{equation}
where
$$
\phi_k\colon\Lambda^k\g^\ndot\otimes M\to M[1-k]
$$
are the Taylor components of the $L_\infty$-module structure.
One easily sees that $T_\pi C_\ndot(A,A)\simeq C_\ndot(A_*,A_*)$
and $T_\pi\Omega^\ndot(\R^d)\simeq\bigl\{\Omega^\ndot(\R^d),d=L_\pi\bigr\}$.

\subsubsection{}
Consider now the disk $D^2$ with the center $\one$ from~[Sh2],
with only one vertex on the boundary, where is placed $\beta\in
A_*$, and $n+1$ points inside, in one of them is placed
$\alpha\in A=T^{-1}_\poly(\R^d)$, and in others are placed
copies of~$\pi$. We can fix the position of $\beta$ because of
the action of the rotation group. Now we consider the
configurations where $\alpha\in[\one,\beta]$, see Figure~1.
There is no edge starting at~$\one$, because we should obtain
a $0$-form.

\sevafig{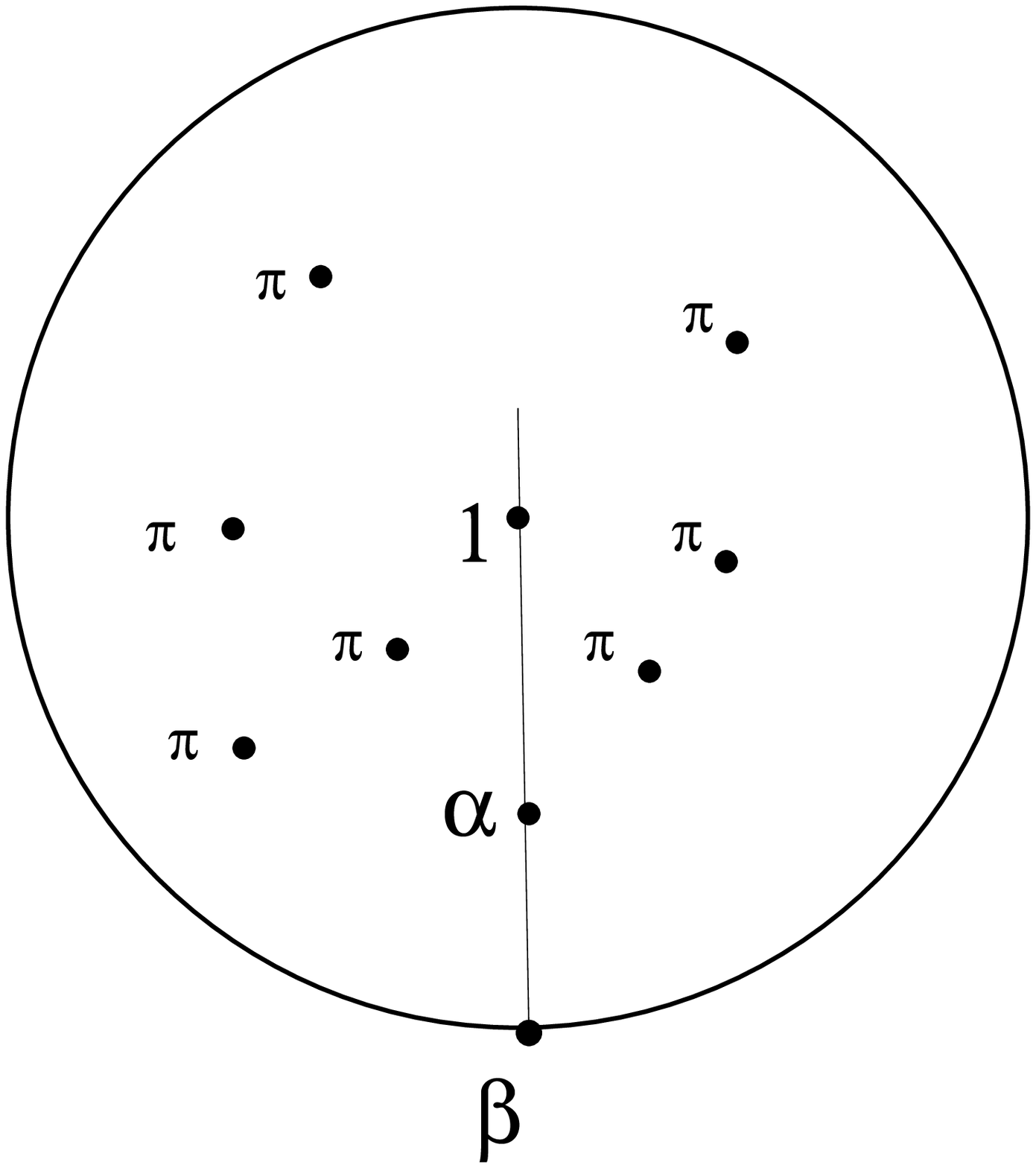}{80mm}{0}{A typical configuration we consider}

We consider the sum over all admissible graphs with
$2(n+1)+1-2=2n+1$ edges, i.e.\ by $1$ less that the usual
configurations in~[Sh2]. But now $\alpha$ moves along the
interval $[\one,\beta]$, and the dimension of the configuration
space is equal to $2n+1$. Denote by $D_{\one,n+1,1}^r$ this
configuration space ($r$~stands for ``restricted''), and
consider any admissible graph $\Gamma$ with $2n$ edges. We have:
\begin{equation}
\int_{\bar D_{\one,n+1,1}^r}d\left(\bigwedge_{e\in E_\Gamma}d\varphi_e\right)=0
\end{equation}
and, by the Stokes formula,
\begin{equation}
\int_{\partial\bar D_{\one,n+1,1}^r}\bigwedge_{e\in E_\Gamma}d\varphi_e=0
\end{equation}
Now we want to describe the boundary strata in $\partial\bar
D_{\one,n+1,1}^r$ of codimension~$1$. There are many
possibilities. First look for the most interesting:

\begin{description}
\item[S1)] the point $\alpha$ and the points
$p_{i_1},\dots,p_{i_k}$ of the first type tend to~$\one$;
\item[S2)] the point $\alpha$ and the points
$p_{i_1},\dots,p_{i_k}$ of the first type tend to
$\beta\in\partial\bar D^2=S^1$;
\item[S3)]  points $p_{i_1},\dots,p_{i_k}$ of the first type,
$p_{i_s}\ne\alpha$ for any~$s$, tend to~$\one$.
\end{description}

There are also other possibilities:

\begin{description}
\item[S4)] points $p_{i_1},\dots,p_{i_k}$ of the first type,
$p_{i_s}\ne\alpha$ for any~$s$, tend to~$\beta$;
\item[S5)] points $p_{i_1},\dots,p_{i_k}$, $p_{i_s}\ne\alpha$
for any~$s$, tend to a point $\kappa$ on the boundary,
$\kappa\ne\beta$;
\item[S6)] points $p_{i_1},\dots,p_{i_k}$, $p_{i_s}\ne\alpha$
for any~$s$, tend to~$\alpha$ and far from~$\one$ and
from~$\beta$;
\item[S7)] points $p_{i_1},\dots,p_{i_k}$, $p_{i_s}\ne\alpha$
for any~$s$, tend to each other inside the disk.
\end{description}

We have:
\begin{equation}
0=\int_{\partial\bar D_{\one,n+1,1}}^r\bigwedge_{e\in E_\Gamma}d\varphi_e=
\int_{\partial_{\text{S1)}}}\natural+
\int_{\partial_{\text{S2)}}}\natural+
\int_{\partial_{\text{S3)}}}\natural+
\int_{\partial_{\text{S4)}}}\natural+
\int_{\partial_{\text{S5)}}}\natural+
\int_{\partial_{\text{S6)}}}\natural+
\int_{\partial_{\text{S7)}}}\natural
\end{equation}
where $\natural=\bigwedge\limits_{e\in E_\Gamma}d\varphi_e$.

We claim, that only
$\int_{\partial_{\text{S1)}}}\natural$,
$\int_{\partial_{\text{S2)}}}\natural$  and $\int_{\partial_{\text{S3)}}}\natural$ are not equal to~$0$, and
$\int_{\partial_{\text{S1)}}}\natural$ gives exactly first summand of the r.h.s.\
of~(6),
$\int_{\partial_{\text{S2)}}}\natural$ gives  the l.h.s.\
of~(6), and $\int_{\partial_{\text{S3)}}}\natural$ gives the second summand in the r.h.s. of (6), $c(\alpha,\beta)$.
Therefore, we consider at first these three cases.

\subsubsection{The case S1)}
By Theorem~6.6.1 in~[K], the integral over this boundary stratum
may not vanish only $k=0$. The situation is like that: only the
point $\alpha$ approaches to the point $\one$ along the interval
connecting $\one$ and~$\beta$. The dimension of this boundary
stratum is equal to~$0$; therefore, there should be no edges
between $\one$ and~$\alpha$. The picture is like in Figure~2.

\sevafig{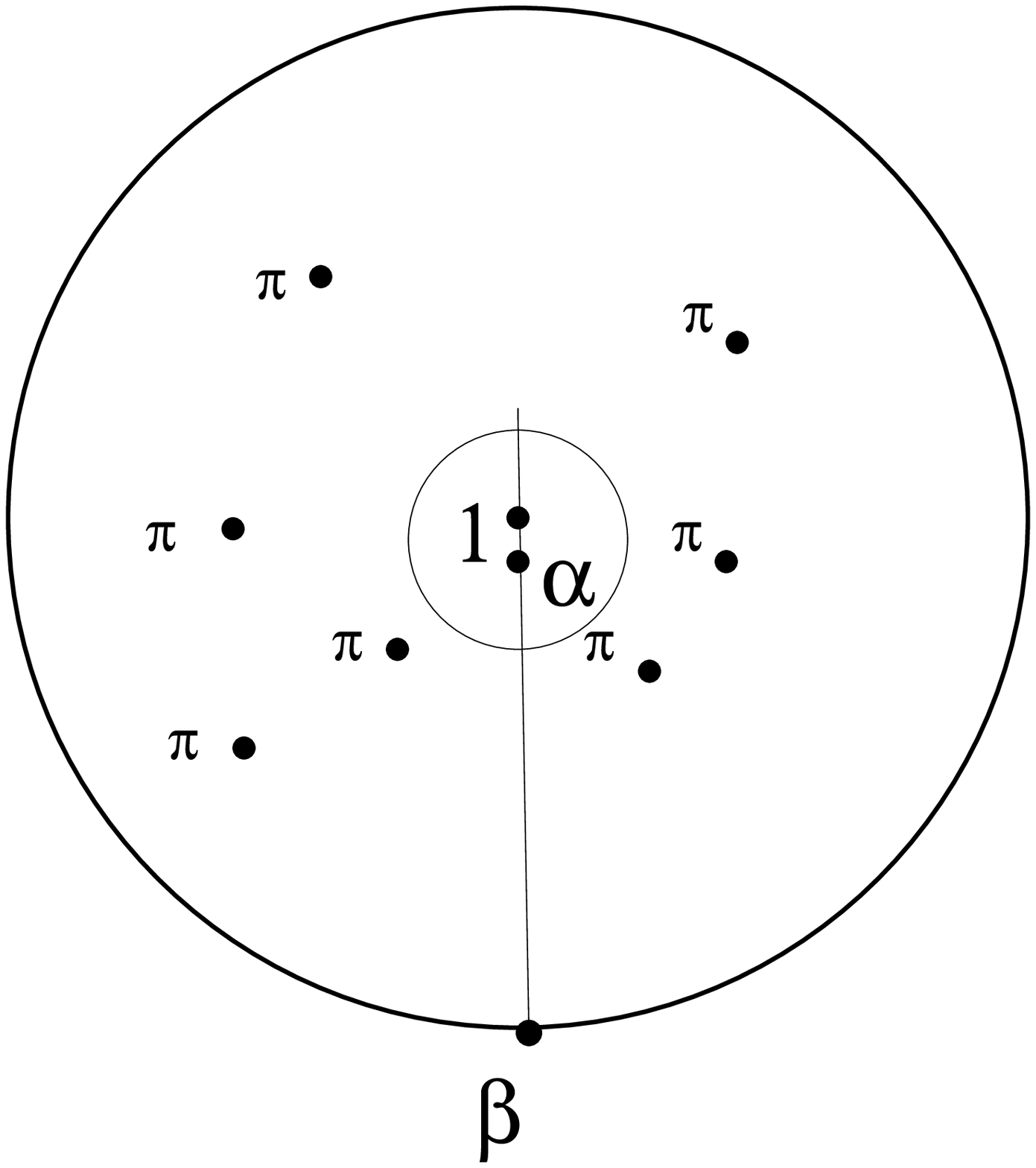}{80mm}{0}{The boundary stratum S1).}

This gives exactly $\alpha\cdot T_\pi\hat\U(\beta)$, i.e. the first summand in the
right-hand side of~(6).

\subsubsection{The case S2)}
In this case $\alpha$ approaches to~$\beta$. Also, some other
points $p_{i_1},\dots,p_{i_k}$ approach to~$\beta$. The
situation can be described in three steps.

\subsubsubsectionb{}
At first, we have the Kontsevich-type picture for this boundary
stratum. It means that we consider the space $C^r_{k+1,1}$
from~[K], where $\alpha$ belongs to a vertical line passing
through~$\beta$ (``$r$''~stands for restricted).

\sevafig{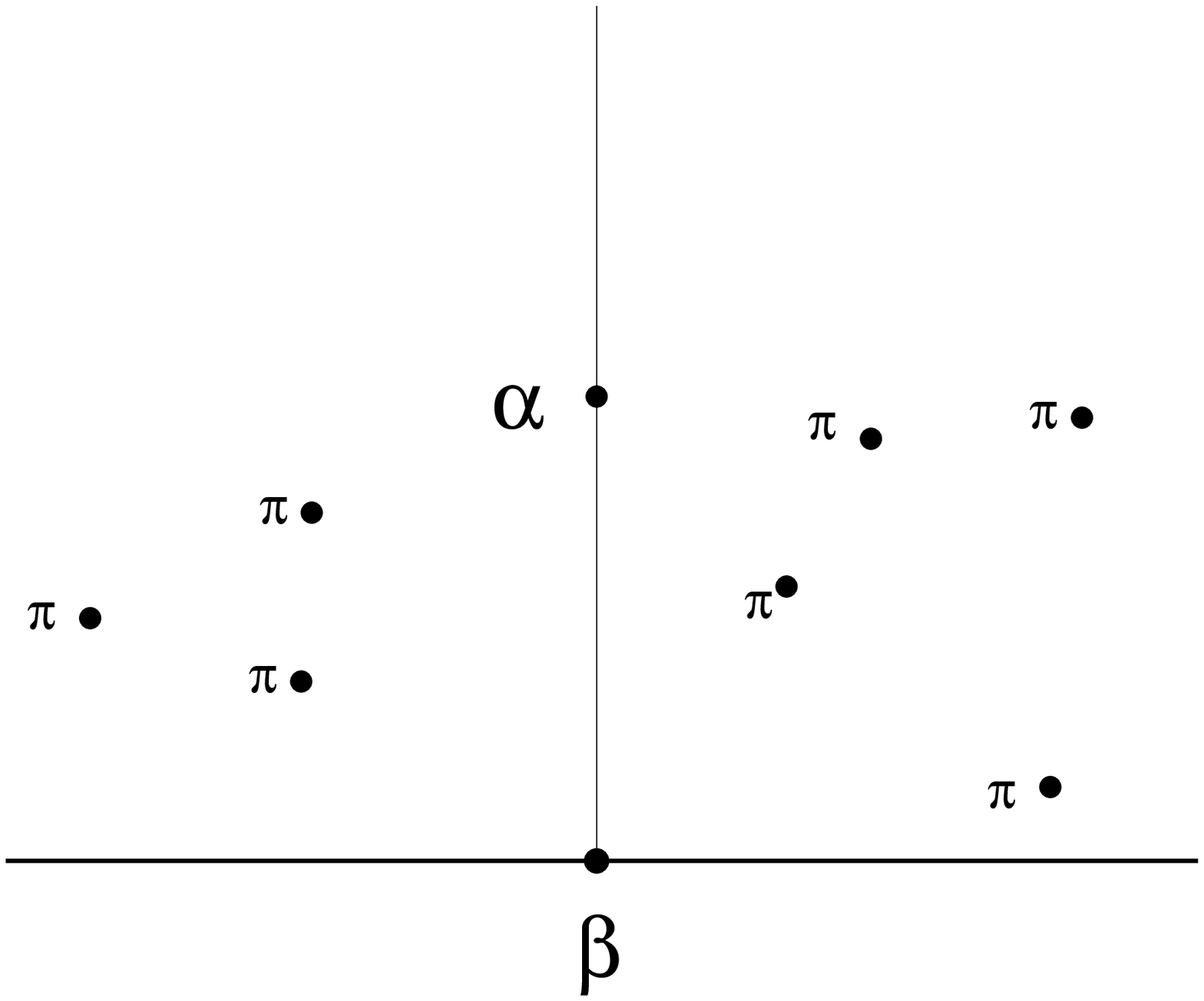}{80mm}{0}{The boundary stratum S2). First reduction}

The dimension of this stratum is by $1$ less than $C_{k+1,1}$,
that is it is equal to $2(k+1)+1-2-1=2k$. Now $\alpha$ is on a
finite distance from the boundary. We want to compute the
corresponding Kontsevich (poly)differential operator. To do
this, we use a second reduction.

\subsubsubsectionb{}
Now we move $\alpha$ to the boundary.

\sevafig{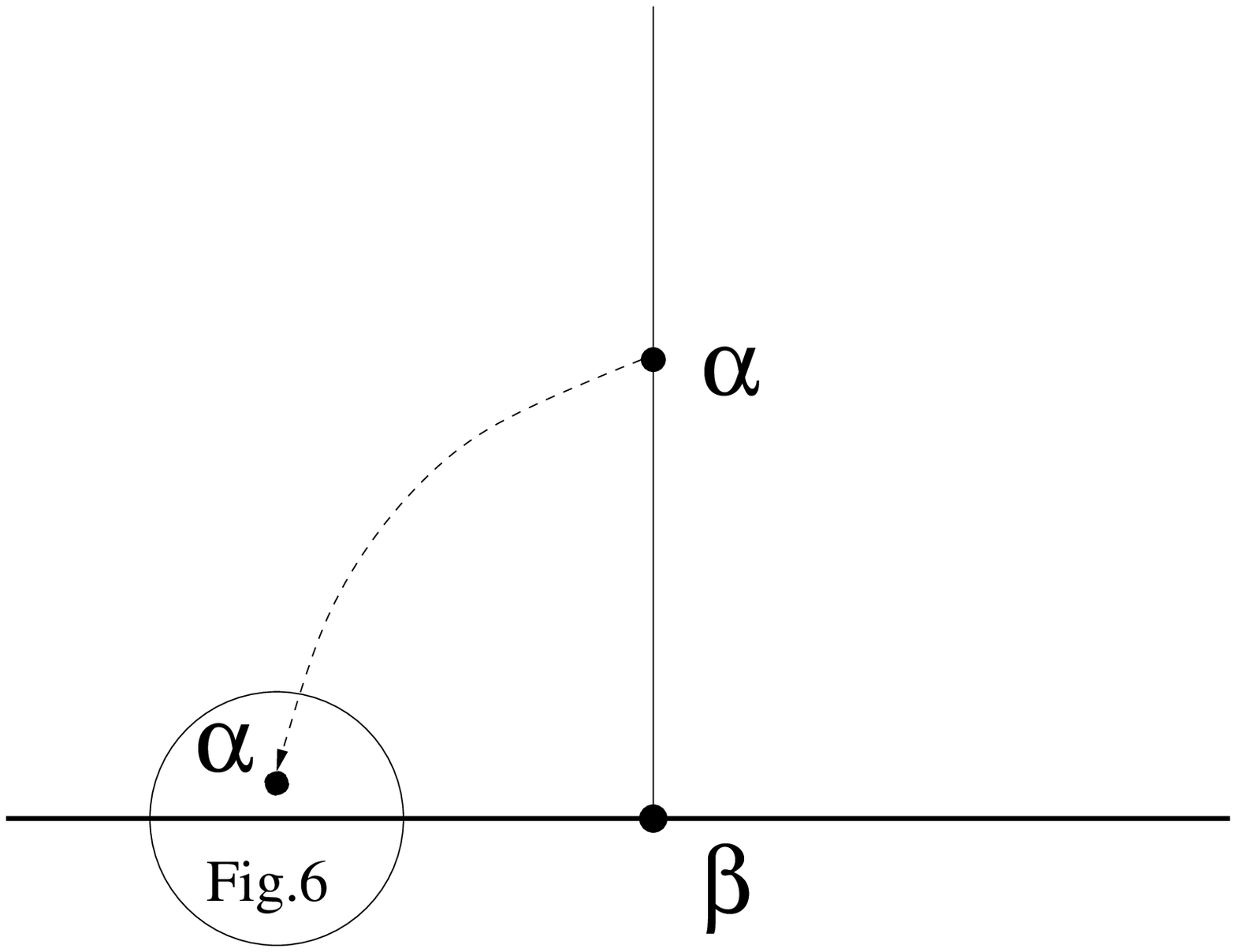}{80mm}{0}{The boundary stratum S2). Second reduction}

In its ``final'' position, $\alpha$ approaches the boundary, but
it still belongs to the interior of the upper half-plane. We
obtain a boundary stratum of codimension~$1$ of $\bar
C_{k+1,1}$, \emph{not}~$\bar C_{k+1,1}$. This boundary stratum
has the same dimension that the space $\bar C^r_{k+1}$, i.e.\
has a codimension~$1$ in $\bar C_{k+1,1}$. We claim, that
polydifferential operators, corresponding to Figure~3 and to
Figure~4, coincide. It's again an application of the Stokes
formula: certainly, we have some other boundary strata, when
several points move close to~$\alpha$ in its intermediate
position, see Figure~5.

\sevafig{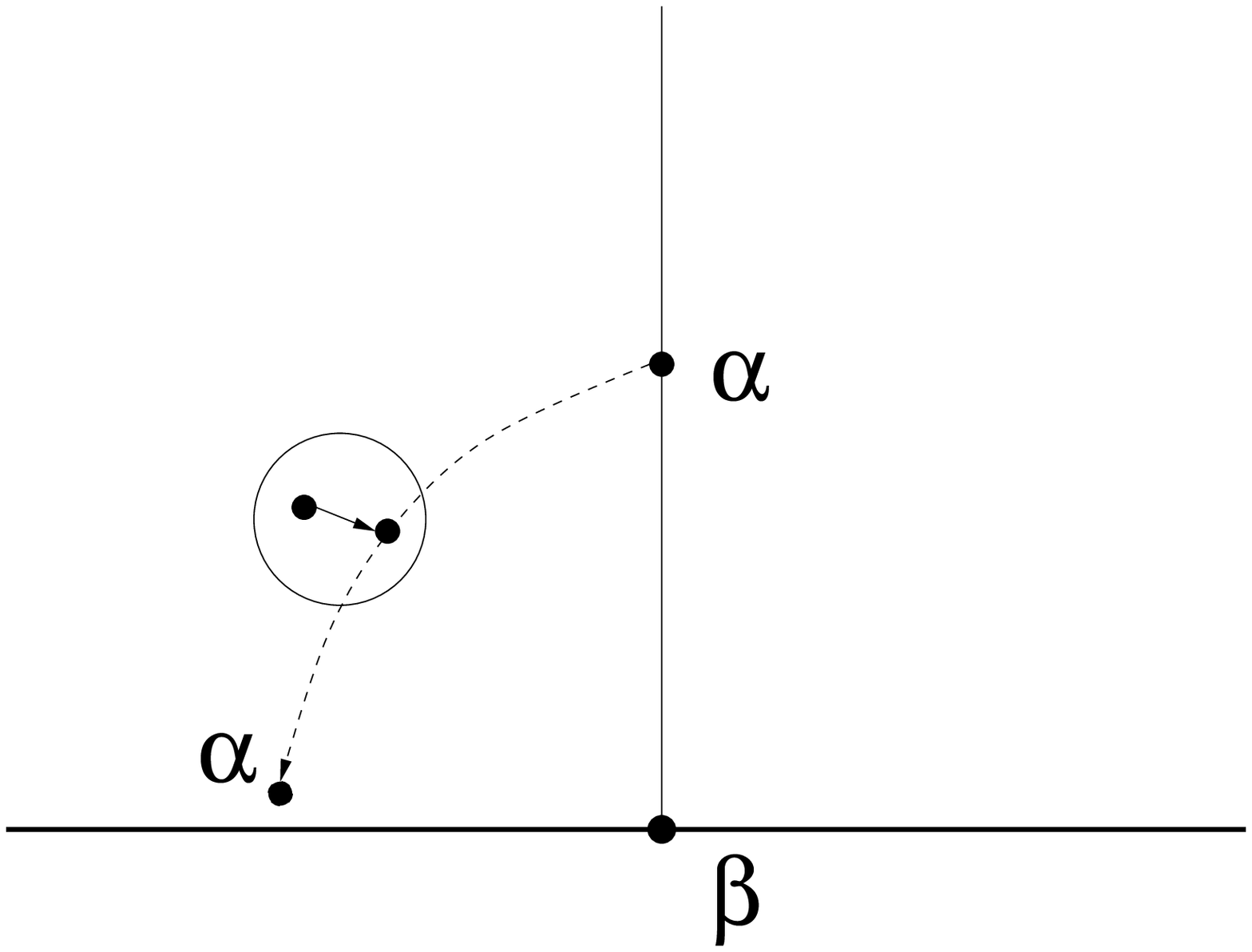}{80mm}{0}{An intermediate position of $\alpha$ and a boundary stratum}

Some points $p_{j_1},\dots,p_{j_l}$ approaches~$\alpha$. By  the
Theorem~6.6.1  from~[K],  we have $l=1$. There should be exactly
one edge from $p_j=p_{j_1}$ to~$\alpha$. This term corresponds
to the bracket $[\pi,\alpha]$, which vanishes because $\alpha$
is supposed to be invariant.

At the picture, showed in Figure~4, we have the polydifferential
operator $T(\alpha,\pi)*\beta$, where $*$ is the Kontsevich
star-product, and $T(\alpha,\pi)$ is an expression, corresponding
to the boundary stratum in Figure~4.

\subsubsubsectionb{}
The picture for $T(\alpha,\pi)$ is showed in Figure~6.

\sevafig{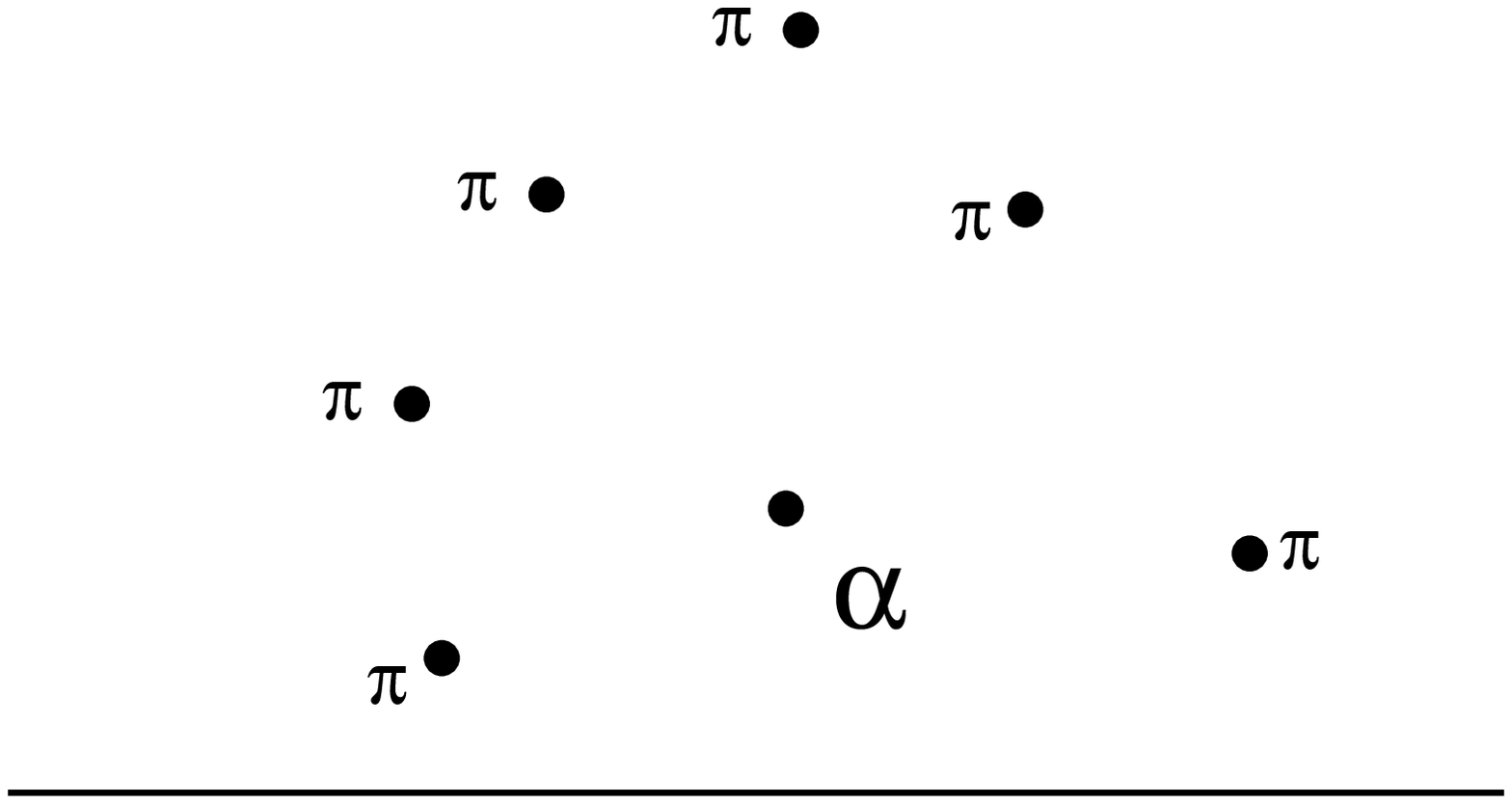}{80mm}{0}{The boundary stratum S2). Third reduction}

It is the usual Kontsevich's picture from~[K]. The corresponding
function $T(\alpha,\pi)$ is equal to~$T_\pi\U(\alpha)$. Finally,
we see that the expression corresponding to the boundary
stratum~S2), is $T_\pi\hat\U((T_\pi\U(\alpha))*\beta)$. The
$T_\pi\hat\U$ outside parentheses is corresponded to Figure~3.

\subsubsubsectionb{}

\begin{remark}
As well we can move $\alpha$ to the right from~$\beta$. We will
obtain $T_\pi\hat\U(\beta*T_\pi\U(\alpha))$. The both
expressions coincide because $\alpha$ satisfies
$[\pi,\alpha]=0$, and, therefore, $T_\pi\U(\alpha)$ is a central
element in the deformed algebra. See~[K], Section~8.
\end{remark}

\subsubsection{The case S3)}

In this case the boundary stratum is $\bar D_k\times\bar
D_{\one,n-k+1,1}^r$, where
$$
D_k=\bigl\{p_1,\dots,p_k\in\C,p_i\ne p_j\ \text{for}\ i\ne
j\bigr\}\big/\bigl\{z\mapsto az,a\in\R_>\bigr\}.
$$
The space $D_k$ has dimension $2k-1$.It follows from Theorem 6.6.1 in [K] that the integral over $D_k$ does not vanish only when $k=1$, and there is no edges from $1$ to $p_1$. (See [Sh2] for some details). This stratum is corresponded to the second summand, $c(\alpha,\beta)$, in the r.h.s. of (6). At the same time, we obtain an explicit formula for $c(\alpha,\beta)$.

\begin{remark}
The stratum S3) here is what was called S2.2) in [Sh2].
\end{remark}

\subsubsection{}

Here we consider the remaining cases S4)-S7).

\subsubsubsection{The case S4)}
In this case the boundary stratum is $\bar C_{k,1}\times\bar
D_{\one,n-k,1}^r$, it has codimension~$1$ as expected. The
integral factories to the product of $a_n$ integral over~$\bar
C_{k,1}$ and an integral over~$\bar D_{\one,n-k,1}^r$. It is
clear that the integral over~$\bar C_{k,1}$ vanishes: we attach
the bivector field~$\pi$ to any point~$p_{i_s}$, therefore, the
number of edges of any graph is~$2k$. But $\dim\bar
C_{k,1}=2k-1$.

\begin{remark}
In the case when $\alpha$ also approaches to~$\beta$ this
argument does not hold, because there are no edges starting
at~$\alpha$, and $\dim\bar C_{k+1,1}^r=2k$.
\end{remark}

\subsubsubsection{The case S5)}
The boundary stratum is $\bar C_{k,0}\times\bar
D_{\one,n-k,2}^r$, it has codimension~$1$. The integral
over~$\bar C_{k,0}$ vanishes because any $p_i$ is a bivector
field, but $\dim\bar C_{k,0}=2k-2$.

\subsubsubsection{The case S6)}
It is the most principal point that this stratum does not
contribute to the integral. By Kontsevich theorem~6.6.1 from~[K]
we have $k=1$. There is only one edge passing from $p_{i_1}$ to
$\alpha$, it corresponds to $[\pi,\alpha]=0$ by the assumption.

\subsubsubsection{The case S7)}
Again, $k=2$ by the Theorem~6.6.1 from~[K]. We have $[\pi,\pi]$
which is equal to~$0$, because $\pi$ is a Poisson bivector
field.

Theorem 2.1 is proven.
\qed

\section{Applications}

For a general Poisson structure $\pi$ on~$\R^d$, the picture is
the following. Denote by $A=C^\infty(\R^d)[[h]]$, by $A_*$ the
Kontsevich deformation quantization of~$A$ (with the harmonic
angle function). Then, we have two maps:
\begin{align*}
T_\pi\U&\colon A\to A_*\\
T_\pi\hat\U&\colon A_*\to A
\end{align*}
such that

\begin{enumerate}
\item
\begin{equation}
(T_\pi\U)(\alpha\cdot\beta)=((T_\pi\U)\alpha)*(T_\pi\U(\beta))
\end{equation}
for $\alpha,\beta$ such that
$
[\pi,\alpha]=[\pi,\beta]=0
$;
in particular, $T_\pi\U$ maps the Poisson center to the center
of the deformed algebra (for the proof see~[K], Section~8);
\item
\begin{equation}
\text{$T_\pi\hat\U$ maps $[A_*,A_*]$ to $\{A,A\}$}
\end{equation}
(for the proof see~[Sh2], Section~3);
\item
compatibility of $T_\pi\U$ and $T_\pi\hat\U$:
\begin{equation}
T_\pi\hat\U(T_\pi\U(\alpha)*\beta)=\alpha\cdot T_\pi\hat\U(\beta)+c(\alpha,\beta)
\end{equation}
here $\alpha\in A$, $[\pi,\alpha]=0$, and $\beta\in A_*$ is
arbitrary, and $c(\alpha,\beta)\in \{A,A\}$ (the same that~(6), is proven in Theorem~2.1 of the
present paper).
\end{enumerate}

Now we are going to consider in more details the case of a
linear Poisson structure.

\subsection{}
Let $\pi$ be a linear Poisson structure on $\R^d\simeq\g^*$,
$\g$~is a finite-dimensional Lie algebra, $\pi$ is the Kostant-Kirillov
Poisson structure. By definition,
$$
\pi=\sum _{ijk=1}^{\dim\g}c_{ij}^k x_k\frac{\partial}{\partial x_i}\wedge
\frac{\partial}{\partial x_j}
$$
where $c_{ij}^k$ are the structure constants of the Lie algebra $\g$ in the
basis $\{x_i\}$.

We proved in [Sh1] that $T_\pi\U=\Id$ in this case, it is useful
(but not necessarily) to use this result here. We have from~(16)
\begin{equation}
T_\pi\hat\U(\alpha*\beta)=\alpha\cdot T_\pi\hat\U(\beta)+c(\alpha,\beta)
\end{equation}
for any $\beta$ and $\alpha$ such that $[\pi,\alpha]=0$. Here
$*$ is the Kontsevich star-product. Now we can suppose that $\g$ is semisimple and we have the decompositions (4) and (5).
Therefore, when we set $\beta=1$ we
obtain
\begin{equation}
T_\pi\hat\U(\alpha)=\alpha
\end{equation}
for any $\alpha$ such that $[\pi,\alpha]=0$. A~priori we have
from~[Sh2]:
\begin{equation}
T_\pi\hat\U(f)=\exp\left(\sum_{k\ge1}w_{2k}\Tr_{2k}\right)(f)
\end{equation}
for some complex numbers $\{w_{2k}\}$. It is enough to know~(17)
for $\alpha$ such that $[\pi,\alpha]=0$, and for $\g=\gl_n$,
$n\ge1$, to conclude that
\begin{equation}
T_\pi\hat\U(\alpha)=\alpha\quad\text{for any}\ \alpha\in A_{(*)}.
\end{equation}
The coefficients $\{w_{2k}\}$ do not depend on the Lie
algebra~$\g$, and we have proved~(19) for any~$\alpha$ and any
Lie algebra~$\g$.

\begin{theorem}
For any finite-dimensional Lie algebra~$\g$, any $\alpha\in
[S^\ndot(\g)]^\g$ and any $\beta\in S^\ndot(\g)$ one has\emph:
\begin{equation}
\alpha*\beta=\alpha\cdot\beta +c(\alpha,\beta)
\end{equation}
where $c(\alpha,\beta)\in \{S^\ndot(\g),S^\ndot(\g)\}$.
\qed
\end{theorem}

\begin{corollary}
For any finite-dimensional Lie algebra~$\g$, any
$\alpha\in[S^\ndot(\g)]^\g$ and any $\beta\in S^\ndot(\g)$ one
has\emph:
\begin{equation}
\varphi_D(\alpha\cdot\beta +c(\alpha,\beta))=\varphi_D(\alpha)*\varphi_D(\beta)
\end{equation}
where $\varphi_D\colon S^\ndot(\g)\to\U(\g)$ is the Duflo map,
and $*$ is the product in~$\U(\g)$.
\end{corollary}

\begin{proof}
The natural isomorphism of \emph{algebras} $\Theta\colon
S(\g)_*\to\U(\g)$,
$$
\Theta(g_1*\dots*g_k)=g_1\otimes\dots\otimes g_k
$$
is equal to $\varphi_D$ (see~[Sh1]). We just apply the map
$\Theta$ to both sides of~(20) and use that $\Theta$ is a map of
algebras.
\end{proof}

\subsection{Remark}
It is an interesting question does the Kashiwara-Vergne
conjecture~[KV] imply our result in Theorem~3.1.
On the other hand, it is interesting does our result (with an  explicit form of
 $c(\alpha,\beta)$) opens a way to prove the Kashiwara-Vergne
conjecture itself.
\section{Acknowledgements}

Giovanni Felder explained to me the proof of the Kontsevich's
theorem on cup-products from~[K], Section~8. Discussions with Boris
Feigin  were very useful for me. I would
like to thank Dominique Manchon for sending me a first version of the paper [MT]. The work was done during my stay at the ETH-Zentrum
in remarkable and stimulating atmosphere. I am grateful to IPDE
grant 1999--2001 for a particular financial support.

 \end{document}